\documentclass[12pt]{article}
\usepackage{a4,bbold}

\newfont{\goth}{eufm10 scaled 1200}
\newcommand{\nc}{\newcommand}
\nc{\ep}{\vspace{0.3cm}}
\nc{\Cl}{\C l}
\nc{\ClE}{\C l(E)}
\nc{\ClF}{\C l(F)}
\nc{\Si}{\Sigma}
\nc{\SE}{\Sigma E}
\nc{\SF}{\Sigma F}
\nc{\SEp}{\Sigma^+ E}
\nc{\SEm}{\Sigma^- E}
\nc{\So}{\Sigma^0}
\nc{\Sl}{\Sigma^1}
\nc{\ga}{\gamma}
\nc{\gaE}{\gamma_E}
\nc{\gaF}{\gamma_F}
\nc{\End}{\mbox{End}}
\nc{\EpF}{E\oplus F}
\nc{\oc}{\omega_\C}
\nc{\Id}{\mbox{Id}}
\nc{\spin}{\mbox{\goth spin}}
\nc{\Spin}{\mbox{Spin}}
\nc{\SO}{\mbox{SO}}
\nc{\Dt}{\tilde{D}}
\nc{\Dh}{\widehat{D}}
\nc{\DM}{D_M^{\Si N}}
\nc{\DMt}{\tilde{D}_M^{\Si N}}
\nc{\vol}{\mbox{vol}}
\nc{\LL}{L^2(M)}
\nc{\LLL}{\mbox{\goth L}}
\nc{\HHH}{\mbox{\goth H}}
\nc{\codim}{\mbox{codim}}
\nc{\grad}{{grad}}

\DeclareSymbolFont{Extrasymb}{U}{msa}{m}{n}
\DeclareMathSymbol\square\mathrel{Extrasymb}{"03}
\def\qed{{\leavevmode\unskip\nobreak\hfil\penalty 50\hskip 1em%
  \hbox{}\nobreak\hfil\lower 1pt\hbox{$\square $\kern-.5pt}\parfillskip 0pt
  \finalhyphendemerits 0\par}}

\newtheorem{theorem}{\bf T{\footnotesize HEOREM}}[section]
\newtheorem{lemma}[theorem]{\bf L{\footnotesize EMMA}}
\newtheorem{corollary}[theorem]{\bf C{\footnotesize OROLLARY}}

\begin{document}

\title{Extrinsic Bounds for Eigenvalues of the\\ Dirac Operator}
\author{Christian B\"ar\thanks{partially supported by SFB 256 and by the
GADGET program of the EU}}
\date{October, 1997}
\maketitle

\begin{abstract}
\noindent
We derive upper eigenvalue bounds for the Dirac operator of a closed
hypersurface in a manifold with Killing spinors such as Euclidean
space, spheres or hyperbolic space.
The bounds involve the Willmore functional.
Relations with the Willmore inequality are briefly discussed.
In higher codimension we obtain bounds on the eigenvalues of the
Dirac operator of the submanifold twisted with the spinor bundle of
the normal bundle.


\end{abstract}

%%%%%%%%%%%%%%%%%%%%%%%%%%%%%%%%%%%%%%%%%%%%%%%%%%%%%%%%%%%%%%%%%%%%%%
%%%%%%%%%%%%%%%%%%%%%%%%%%%%%%%%%%%%%%%%%%%%%%%%%%%%%%%%%%%%%%%%%%%%%%

\setcounter{section}{-1}

\section{Introduction}

Lower and upper eigenvalue estimates for operators like the Dirac
operator on a closed Riemannian spin manifold are derived
by very different methods.
Lower estimates are usually based on a Bochner-Weitzenb\"ock formula,
i.e.\ on an intelligent partial integration.
The first idea to obtain upper bounds on Dirac eigenvalues is due to
Vafa and Witten \cite{VW,Atiyah} who show that there are upper
eigenvalue estimates for all twisted Dirac operators on a closed
Riemannian spin manifold solely in geometric data of the manifold
but independent of the twist.

The idea is as follows.
Compare the Dirac operator $D$ (or a multiple of it) to a
twisted Dirac operator $\Dt$ acting on sections of the same vector
bundle.
By index theory make sure that $\Dt$ has a kernel.
Let $k$ be the multiplicity of the eigenvalue 0 of $\Dt$.
Estimate the norm of the difference (which is of zero order),
$\|D-\Dt\|=:C$, by geometric quantities.
Then at least $k$ eigenvalues of $D$ are bounded by $C$.

How can one get a good twisted Dirac operator $\Dt$ to compare $D$
with?
One way to achieve this is to take a map $f:M\to S^n$ of nonzero
degree and pull-back suitable vector bundles from the sphere to
twist the Dirac operator with.
The fact that $S^n$ has no cohomology in middle dimensions helps a lot
for the index computations.
Taking for $f$ the Gauss map of a surface composed with a suitable
self-mapping of $S^2$ Baum used this idea to show
\ep
\begin{theorem} \label{baumtheo}
{\em (Baum \cite[Prop.~2]{Baum})}
Let $M\subset\R^3$ be a compact oriented surface of genus $g\not= 1$.
Then there is an eigenvalue $\lambda$ of the Dirac operator on $M$
satisfying
$$
|\lambda| \le c(g) \cdot \max\{\mbox{principal curvatures of }M\}
$$
where
$$
c(g) = \left\{
\begin{array}{ll}
1 & \mbox{if }g=0 \\
3 & \mbox{if }g=2,3 \\
2 & \mbox{if }g\ge 4
\end{array}
\right.
$$
\end{theorem}
\ep

The torus ($g=1$) cannot be dealt with by this approach because its
Gauss map has degree zero.
Baum could also use this method to obtain intrinsic upper bounds
on the Dirac eigenvalues of a compact Riemannian spin manifold
in terms of sectional curvature and injectivity radius
\cite[Prop.~1]{Baum}.

A variation of the same approach was used by Bunke to show
\ep

\begin{theorem} \label{bunketheo}
{\em (Bunke \cite[Thm.~A]{Bunke})}
Let $M$ be an $n$-dimensional compact Riemannian spin manifold
isometrically immersed in $\R^N$.
Let $II$ denote its second fundamental form.

Then there is a topologically determined number of Dirac eigenvalues
of $M$ satisfying
$$
\lambda^2 \le 2^{[n/2]}\cdot \|II\|^2_{L^\infty(M)} .
$$
\end{theorem}
\ep

The topologically determined number of Dirac eigenvalues which can be
estimated can be given explicitly in terms of $N$ and of characteristic
numbers of $M$.
It can be zero however.
So the theorem is not always applicable.
For example, the case of a torus in $\R^3$ cannot be handled by
this theorem either.

The other main method to derive upper eigenvalue bounds is based on the
variational characterization of eigenvalues.
If one has a $k$-dimensional space of ``test spinors'' $\phi$ on
which the {\em Rayleigh quotient} $(D^2\phi,\phi)_{L^2}/(\phi,\phi)_{L^2}$
is bounded by some constant $C$, then there are at least $k$
eigenvalues of $D^2$ bounded by $C$.

This approach has been used by the author in \cite{Baer} to get
intrinsic upper bounds on Dirac eigenvalues in terms of sectional
curvature and injectivity radius.
These estimates are sharp for spheres of constant curvature.

Anghel obtained the following bound on {\em spectral gaps}.
If we order the Dirac eigenvalues by increasing absolute value,
$0 \le |\lambda_1| \le |\lambda_2| \le \cdots \to \infty$, then
\ep

\begin{theorem}
{\em (Anghel \cite[Thm.~3.1]{Anghel})}
Let $M$ be an $n$-dimensional compact Riemannian spin manifold
isometrically immersed in $\R^N$.
Let $H$ denote its mean curvature vector field.
Let $S_0$ denote the minimum of the scalar curvature of $M$.

Then one has
$$
\lambda_{m+1}^2 - \lambda_m^2 \le
n\|H\|^2_{L^\infty(M)} + \frac{4}{mn}\sum_{k=1}^m \lambda_k^2
-\frac{S_0}{n}.
$$
\end{theorem}
\ep

The proof is based on the variational characterization of eigenvalues.
One constructs test spinors for the $(m+1)^{st}$ eigenvalue using the
eigenspinors of the previous eigenvalues and the coordinate functions
given by the immersion.

To convert this into an upper bound on the eigenvalues themselves
one has to assume something on the smallest eigenvalue.
For example, if 0 is an eigenvalue and if the scalar curvature
vanishes identically, $S\equiv 0$, then one concludes for the smallest
nonzero eigenvalue $\lambda$
\cite[Thm.~3.6]{Anghel}
$$
\lambda^2 \le \frac{n}{\vol(M)}\int_M |H|^2.
$$

We will show that on a hypersurface of $\R^{n+1}$ a certain number of
Dirac eigenvalues can always be bounded in terms of $\int_M H^2$
without any a-priori assumption on the spectrum or on the scalar
curvature.
More precisely, we will prove
\ep

\noindent
{\bf C{\footnotesize OROLLARY} \ref{hypeukl}}
{\em Let $M$ be an $n$-dimensional closed oriented hypersurface
isometrically immersed in $\R^{n+1}$.
Let $M$ carry the induced spin structure.
Let $H$ be the mean curvature of $M$ in $\R^{n+1}$.

Then the classical Dirac operator $D_M$ of $M$
has at least $2^{[\frac{n}{2}]}$ eigenvalues
$\lambda$ (counted with multiplicities) satisfying}
$$
\lambda^2 \le \frac{n^2}{4\cdot\vol(M)} \int_M H^2 .
%\mbox{\qed}
$$
\ep

Note that this improves Theorem \ref{baumtheo} since the maximum
of the principal curvatures is replaced by $\int_M H^2/\vol(M)$, the
constant $c(g)$ is replaced by 1, and the torus case is included.

We can replace the ambient space $\R^{n+1}$ by other spaces like
the standard sphere $S^{n+1}$ or hyperbolic space $H^{n+1}$ in
which cases we obtain
\ep

\noindent
{\bf C{\footnotesize OROLLARY} \ref{hypsphaere}}
{\em Let $M$ be an $n$-dimensional closed oriented hypersurface
isometrically immersed in $S^{n+1}$.
Let $M$ carry the induced spin structure.
Let $H$ denote the mean curvature of $M$ in $S^{n+1}$.

Then the classical Dirac operator $D_M$ of $M$
has at least $2^{[\frac{n}{2}]}$ eigenvalues
$\lambda$ satisfying}
$$
\lambda^2 \le \frac{n^2}{4} + \frac{n^2}{4\cdot\vol(M)} \int_M H^2
%\mbox{\qed}
$$
\ep

\noindent
and
\ep

\noindent
{\bf C{\footnotesize OROLLARY} \ref{hyphyp}}
{\em Let $M$ be an $n$-dimensional closed oriented hypersurface
isometrically immersed in $H^{n+1}$.
Let $M$ carry the induced spin structure.
Let $H$ denote the mean curvature of $M$ in $H^{n+1}$.

Then the classical Dirac operator $D_M$ of $M$
has at least $2^{[\frac{n}{2}]}$ eigenvalues
$\lambda$ satisfying}
$$
|\lambda| \le \frac{n}{2}\left( 1 + \| H\|_{L^\infty(M)} \right).
%\mbox{\qed}
$$
\ep

\noindent
Equality is attained in corollaries \ref{hypeukl} and \ref{hypsphaere}
for all distance spheres.

In general, the ambient space can be any Riemannian spin manifold
carrying Killing spinors.
These Killing spinors are then restricted to the hypersurface and
used as test spinors.
In case the hypersurface bounds an open subset in the ambient space
more eigenvalues can be bounded (Theorem \ref{hoehere}).
This is sharp and optimal as the example of the standard sphere
in $\R^{n+1}$ shows.
The statement of Theorem \ref{hoehere} is false if the hypersurface
does not bound.

It is also possible to look at submanifolds of higher codimension.
Then we obtain upper bounds on the eigenvalues of the Dirac operator
on the submanifold twisted with the spinor bundle of the normal
bundle.
In case of a hypersurface our estimates improve those of Theorem
\ref{bunketheo}, in case of higher codimension they are logically
independent because they make statements about different operators.

The paper is organized as follows.
In the first section we study spinor modules for the Clifford
algebra of a direct sum of two Euclidean vector spaces.
Later this will be applied to the sum of the tangent space and
the normal space of the submanifold.

In the second part we compare the spinor connection of the submanifold
to the spinor connection of the ambient space.
This implies a relation between the Dirac operator of the ambient
space and the Dirac operator of the submanifold twisted with the
spinor bundle of the normal bundle.
The mean curvature appears as a correction term.
We hope that the first two sections in which the ``submanifold
theory'' of Dirac operators is established will also be of
independent interest.

In the third section we prove the eigenvalue estimate for submanifolds
of arbitrary codimension in a Riemannian spin manifold with
Killing spinors.
The main result for real Killing constant is Theorem \ref{abschreell}.
The case of imaginary Killing constant is somewhat more complicated
because then Killing spinors do not have constant length.
The estimate is given in Theorem \ref{abschimaginaer}.

In the forth part we restrict our attention to hypersurfaces.
The results of the previous section yield bounds on the
eigenvalues of the classical (untwisted) Dirac operator of the
hypersurface.
We also discuss relations with the Willmore problem in surface theory.
Moreover, we explain how one can bound more eigenvalues if (and only
if) the hypersurface bounds an open subset.

In the last section we show how to get upper bounds on {\em all} Dirac
eigenvalues.
These estimates involve the Laplace eigenvalues of the hypersurface.
\ep

{\bf Acknowledgement.}
It is a pleasure to thank Bernd Ammann for helpful discussions.

%%%%%%%%%%%%%%%%%%%%%%%%%%%%%%%%%%%%%%%%%%%%%%%%%%%%%%%%%%%%%%%%%%%%%%
%%%%%%%%%%%%%%%%%%%%%%%%%%%%%%%%%%%%%%%%%%%%%%%%%%%%%%%%%%%%%%%%%%%%%%

\section{Algebraic Preliminaries}

Our aim is to compare the Dirac operator on a Riemannian spin manifold
with the one on a spin submanifold.
In particular, we have to compare the restriction of the spinor bundle
to the submanifold with the genuine spinor bundle of the submanifold
itself.
The starting point is the splitting of the restricted tangent bundle
of the large manifold into tangent and normal bundle of the submanifold.
Hence we need to compare the spinor modules of the Clifford algebra
of a direct sum of two Euclidean vector spaces with the spinor modules
associated with the two factors.
In principle, this is contained in \cite[Ch.~I.5]{LM} but for our
considerations we have to make things more explicit.

Let $E$ be an oriented Euclidean vector space.
We denote the complex Clifford algebra of $E$ by $\ClE$.
For basics on Clifford algebras and spinors see e.g.\ \cite{BGV} or
\cite{LM}.
If the dimension $n$ of $E$ is even, then $\ClE$ has precisely one
irreducible module, the spinor module $\SE$.
It has $\dim(\SE)=2^{n/2}$.
Denote the Clifford multiplication by $\gaE : \ClE \to \End(\SE)$.
When restricted to the even subalgebra $\C l^0(E)$ the spinor
module decomposes into even and odd half-spinors $\SE = \SEp \oplus \SEm$.
The ``complex volume element'' $\oc = i^{n/2} \gaE(e_1\cdots e_n)$
acts as $+1$ on $\SEp$ and as $-1$ on $\SEm$.
Here $e_1,\ldots,e_n$ denote a positively oriented orthonormal basis of
$E$.

If $n$ is odd there are exactly two irreducible
modules, $\So E$ and $\Sl E$, again called spinor modules.
In this case $\dim(\Sigma^{0,1}E)=2^{(n-1)/2}$.
Clifford multiplication will now be denoted by
$\ga_{E,j} : \ClE \to \End(\Sigma^{j}E)$.
Similarly to the half-spinor spaces in even dimensions the two modules
$\So E$ and $\Sl E$ can be distinguished by the action
of the complex volume element $\oc = i^{(n+1)/2} \gaE(e_1\cdots e_n)$.
On $\Si^jE$ it acts as $(-1)^j$, $j=0,1$.

One can pass from $\So E$ to $\Sl E$ by taking the same underlying vector
space, $\So E = \Sl E$, and setting $\ga_{E,1}(X) := -\ga_{E,0}(X)$
for all $X\in E$.
In other words, there exists a vector space isomorphism
$\Phi : \Si^0E \to \Si^1E$ such that $\Phi\circ\ga_{E,0}(X) = -\ga_{E,1}(X)
\circ\Phi$ for all $X\in E$.

Now let $E$ and $F$ be two oriented Euclidean vector spaces.
Let $\dim E=n$ and $\dim F=m$.
We want to construct the spinor module(s) of $E \oplus F$ from those
of $E$ and $F$.
\ep

\newpage

{\bf Case 1.}
$n$ and $m$ are even.

\noindent
To obtain a vector space of the correct dimension we simply put
$$
\Sigma := \SE \otimes \SF
$$
and
$$
\ga : E \oplus F \rightarrow \End(\Si)
$$
\begin{eqnarray}
\ga(X)(\sigma\otimes\tau) &:=& (\gaE(X)\sigma)\otimes\tau ,
\label{case1} \\
\ga(Y)(\sigma\otimes\tau) &:=& (-1)^{\deg\sigma}\sigma\otimes\gaF(Y)\tau,
\nonumber
\end{eqnarray}
where $X\in E, Y\in F, \sigma\in \SE, \tau\in\SF$ and
$$
\deg\sigma = \left\{
\begin{array}{rl}
0, & \sigma\in\SEp\\
1, & \sigma\in\SEm
\end{array}
\right. .
$$
Here the degree $\deg\sigma$ is such that
$\gamma_E(\oc)\sigma=(-1)^{\deg\sigma}\sigma$.
One easily checks
$$
\ga(X+Y)\ga(X+Y)(\sigma\otimes\tau) = - |X+Y|^2\sigma\otimes\tau.
$$
Thus $\ga$ extends to a homomorphism $\ga : \C l(E\oplus F) \to \End(\Si)$.
Therefore $(\Si,\ga)$ is a nontrivial $\C l(E\oplus F)$-module
of dimension $2^{n/2}\cdot2^{m/2} = 2^{(n+m)/2}$, hence isomorphic to
$(\Si(E\oplus F),\ga_{E\oplus F})$.
The splitting into half-spinor modules is given by
\begin{eqnarray*}
\Si^+(\EpF) &=& (\Si^+E \otimes \Si^+E) \oplus (\Si^-E \otimes \Si^-E)\\
\Si^-(\EpF) &=& (\Si^+E \otimes \Si^-E) \oplus (\Si^+E \otimes \Si^-E).
\end{eqnarray*}
\ep

{\bf Case 2.}
$n$ is even and $m$ is odd.

\noindent
We put
$$
\Si^j := \SE \otimes \Si^jF
$$
for $j=0,1$.
With the same definition of
$$
\ga_j : E \oplus F \rightarrow \End(\Si^{j})
$$
as in Case 1 we again make $\Si^{0}$ and $\Si^{1}$ into
$\C l(\EpF)$-modules.
One easily checks that the complex volume element of $\C l(\EpF)$
acts on $\Si^j$ as $(-1)^j$.
Hence $(\Si^j,\ga_j)$ is isomorphic to $(\Si^j(\EpF),\ga_{\EpF,j})$.
\ep

{\bf Case 3.}
$n$ is odd and $m$ is even.

\noindent
Of course, this case is symmetric to the second one.
Later on we will apply these preliminary considerations to $E$ the tangent
space of a submanifold and to $F$ its normal space.
Then $E$ and $F$ cannot be interchanged.
Therefore let us briefly give the explicit formulas in this case too.

We put
$$
\Si^j := \Si^jE \otimes \SF
$$
and
\begin{eqnarray}
\ga_j : \EpF &\rightarrow& \End(\Si^j) \nonumber \\
\ga_j(X)(\sigma\otimes\tau) &:=& (-1)^{\deg\tau}\cdot\ga_E(X)\sigma
\otimes\tau
\label{case3} \\
\ga_j(Y)(\sigma\otimes\tau) &:=& \sigma\otimes\ga_F(Y)\tau, \nonumber
\end{eqnarray}
where $X\in E$ and $Y\in F$.
One checks that $(\Si^j,\ga_j)$ is a realization of $(\Si^{j}(\EpF),
\ga_{\EpF,j})$.
\ep

{\bf Case 4.}
$n$ and $m$ are odd.

\noindent
This is the most complicated case.
We set
\begin{eqnarray*}
\Si^+ &:=& \Si^0E\otimes\Si^0F,\\
\Si^- &:=& \Si^0E\otimes\Si^1F,\\
\Si&:=&\Si^+\oplus\Si^-.
\end{eqnarray*}
Recall that there exists a vector space isomorphism
$\Phi : \Si^0F \to \Si^1F$ such that $\Phi\circ\ga_{F,0}(Y) = -\ga_{F,1}(Y)
\circ\Phi$ for all $Y\in F$.
With respect to the splitting $\Si = \Si^+\oplus\Si^-$ define
$$
\ga : \EpF \rightarrow \End(\Si), \\
$$
\begin{eqnarray}
\ga(X) &:=& i\cdot
\left(
\begin{array}{cc}
0 & \ga_{E,0}(X)\otimes\Phi^{-1} \\
-\ga_{E,0}(X)\otimes\Phi & 0
\end{array}
\right),
\label{case4} \\
\ga(Y) &:=&
\left(
\begin{array}{cc}
0 & \Id\otimes(\Phi^{-1}\circ\ga_{F,1}(Y)) \\
-\Id\otimes(\Phi\circ\ga_{F,0}(Y)) & 0
\end{array}
\right). \nonumber
\end{eqnarray}
One computes
$$
\ga(X+Y)\circ\ga(X+Y) = -|X+Y|^2\cdot\Id
$$
showing that $\ga$ extends to a representation of $\C l(\EpF)$ on $\Si$.
Moreover, the complex volume element of $\C l(\EpF)$ acts on $\Si^+$ as
$+1$ and on $\Si^-$ as $-1$.

Therefore $(\Si,\ga)$ is a realization of $(\Si(\EpF),\ga_{\EpF})$
and $\Si = \Si^+\oplus\Si^-$ is the splitting into half-spinor spaces.

%%%%%%%%%%%%%%%%%%%%%%%%%%%%%%%%%%%%%%%%%%%%%%%%%%%%%%%%%%%%%%%%%%%%%%
%%%%%%%%%%%%%%%%%%%%%%%%%%%%%%%%%%%%%%%%%%%%%%%%%%%%%%%%%%%%%%%%%%%%%%

\section{Levi-Civita Connection and Dirac Operator}

Let $Q$ be an $(n+m)$-dimensional Riemannian manifold and let $M
\hookrightarrow Q$ be an $n$-dimensional immersed submanifold.
Let $M$ carry the induced Riemannian metric.
We suppose that both manifolds are equipped with a spin structure.
This induces a unique spin structure on the normal bundle $N$ of $M$ in $Q$
such that the sum of the spin structures on the tangent bundle and on the
normal bundle of $M$ coincides with the spin structure
on the tangent bundle of $Q$ restricted to $M$ \cite{Milnor}.
Note that in particular $M$ and $Q$ are oriented.

Denote the Levi-Civita connections of $M$ and $Q$ by $\nabla^M$ and $\nabla^Q$
resp.
Let $\nabla^N$ be the normal connection on $N$.
The second fundamental form of $M$ in $Q$ is denoted by $II$.
For $p \in M$ and $X \in T_pM$ the {\em Gauss formula} says,
with respect to the splitting $T_pQ = T_pM \oplus N_p$,
\begin{equation}
\nabla_X^Q =
\left(
\begin{array}{cc}
\nabla^M_X & -II(X,\cdot)^\ast \\
II(X,\cdot) & \nabla^N_X
\end{array}
\right) .
\end{equation}
In other words,
\begin{equation}
\nabla_X^Q - (\nabla^M_X \oplus \nabla^N_X) =
\left(
\begin{array}{cc}
0 & -II(X,\cdot)^\ast \\
II(X,\cdot) & 0
\end{array}
\right) .
\label{gaussformel2}
\end{equation}

Let $X_1,\ldots,X_n$ be a local positively oriented orthonormal tangent
frame of $M$ near $p$, let $Y_1,\ldots,Y_m$ be a local positively
oriented orthonormal frame of $N$ near $p$.
Then $h := (X_1,\ldots,X_n,Y_1,\ldots,Y_m)$ is a local section of the
frame bundle of $Q$ restricted to $M$.
Now we can write (\ref{gaussformel2}) in matrix form as
\begin{equation}
\nabla_X^Q - (\nabla^M_X \oplus \nabla^N_X) =
\left(
\begin{array}{cc}
0 & (-\langle II(X,X_i),Y_j\rangle)_{j,i} \\
(\langle II(X,X_i),Y_j\rangle)_{i,j} & 0
\end{array}
\right)
\label{gaussformel3}
\end{equation}
Let $\omega^M$, $\omega^N$, and $\omega^Q$ be the connection 1-forms
of $\nabla^M$, $\nabla^N$, and $\nabla^Q$ lifted to $\spin(n)$,
$\spin(m)$, and $\spin(m+n)$ resp.
If $\Theta : \Spin(n+m) \to \SO(n+m)$ is the usual double covering map,
then (\ref{gaussformel3}) can be rewritten as
$$
\Theta_\ast
\left(\omega^Q(dh\cdot X) - (\omega^M \oplus \omega^N)(dh\cdot X)\right) =
$$
\begin{equation}
\left(
\begin{array}{cc}
0 & (-\langle II(X,X_i),Y_j\rangle)_{j,i} \\
(\langle II(X,X_i),Y_j\rangle)_{i,j} & 0
\end{array}
\right) .
\label{gaussformel4}
\end{equation}
Using a standard formula \cite[p.~42]{LM} for $\Theta_\ast$,
(\ref{gaussformel4}) yields
\begin{equation}
\omega^Q(dh\cdot X) - (\omega^M \oplus \omega^N)(dh\cdot X) =
\frac{1}{2} \sum_{i=1}^n \sum_{j=1}^m \langle II(X,X_i),Y_j\rangle
\cdot e_i\cdot f_j
\label{gaussformel5}
\end{equation}
where $e_1,\ldots,e_n$ is the standard basis of $\R^n$,
$f_1,\ldots,f_m$ is the standard basis of $\R^m$ and the
whole expression in (\ref{gaussformel5}) is an element of
$\spin(n+m) \subset \Cl(n+m)$.

>From the considerations in the previous section we know for the spinor
bundles that $\Si Q|_M = \Si M \otimes \Si N$ unless $n$ and $m$ are
both odd in which case $\Si Q|_M = (\Si M \otimes \Si N) \oplus
(\Si M \otimes \Si N)$.

Let $\nabla^{\Si Q}$, $\nabla^{\Si M}$, and $\nabla^{\Si N}$ be the
Levi-Civita connections on $\Si Q$, on $\Si M$, and on $\Si N$.
By
$$
\nabla^{\Si M \otimes \Si N} := \nabla^{\Si M} \otimes \Id +
\Id \otimes \nabla^{\Si N}
$$
we mean the product connection on $\Si M \otimes \Si N$ and also on
$(\Si M \otimes \Si N) \oplus (\Si M \otimes \Si N)$ if $n$ and $m$
are both odd.
Equation (\ref{gaussformel5}) yields
\begin{eqnarray}
\nabla^{\Si Q}_X - (\nabla^{\Si M}_X \otimes \Id + \Id \otimes
\nabla^{\Si N}_X) &=&
\frac{1}{2} \sum_{i=1}^n \sum_{j=1}^m \langle II(X,X_i),Y_j\rangle
\gamma_Q(X_i\cdot Y_j)
\nonumber\\
&=&
\frac{1}{2} \sum_{i=1}^n \gamma_Q(X_i\cdot II(X,X_i)) .
\label{gaussformel6}
\end{eqnarray}

Our aim is to relate the Dirac operator on $M$ twisted with the spinor
bundle of the normal bundle, $D_M^{\Si N}$, to extrinsic quantities on $Q$.
To this extent we first look at the operators
$$
\Dt := \sum_{j=1}^n \gamma_Q(X_j)\nabla_{X_j}^{\Si M\otimes \Si N}
$$
and
$$
\Dh := \sum_{j=1}^n \gamma_Q(X_j)\nabla_{X_j}^{\Si Q} .
$$
It is easy to see that the definitions are independent of the choice
of orthonormal frame $X_1,\ldots,X_n$.
Both operators act on sections of $\Si Q|_M$.
In contrast to the Dirac operator on $M$ they both use Clifford
multiplication $\gamma_Q$ instead of $\gamma_M$.

By $\DM$ we mean the Dirac operator on $M$ twisted with the bundle $\Si N$.
We define
$$
\DMt := \left\{
\begin{array}{cl}
\DM \oplus -\DM & , \mbox{ if $n$ and $m$ are odd} \\
\DM             & , \mbox{ otherwise}
\end{array}
\right.
$$
With this definition $\DM$ acts on sections of $\Si M \otimes \Si N$
and $\DMt$ acts on sections of $\Si Q|_M$.

Let $H = \frac{1}{n}\sum_{j=1}^n II(X_j,X_j)$ be the {\em mean
curvature vector field} of $M$.
\ep

\begin{lemma}\label{diraclemma}
The operator $\Dt$ is formally self-adjoint and we have
$$
\Dt^2 = (\DMt)^2 .
$$
Moreover,
$$
\Dt = \Dh + \frac{n}{2} \gamma_Q(H) .
$$
\end{lemma}
\ep

{\bf Proof.}
Using (\ref{gaussformel6}) we compute
\begin{eqnarray*}
\Dh - \Dt &=& \frac{1}{2} \sum_{i,j=1}^n \gamma_Q(X_j)\gamma_Q(X_i
\cdot II(X_j,X_i)) \\
 &=& \frac{1}{2} \sum_{i,j=1}^n \gamma_Q(X_j\cdot X_i)\gamma_Q(II(X_j,X_i)) .
\end{eqnarray*}
The terms with $i\not= j$ cancel because $X_j\cdot X_i$ is antisymmetric
in $j$ and $i$ whereas $II(X_j,X_i)$ is symmetric.
>From $\gamma_Q(X_j\cdot X_j) = -1$ we obtain
\begin{eqnarray*}
\Dh - \Dt &=& - \frac{1}{2} \sum_{j=1}^n \gamma_Q(II(X_j,X_j)) \\
&=& - \frac{n}{2} \gamma_Q(H) .
\end{eqnarray*}
Hence we have shown
$$
\Dt = \Dh + \frac{n}{2} \gamma_Q(H) .
$$
To relate $\Dt$ and $\DMt$ we have to distinguish the various cases
of section 1.
If $n$ is even (case 1 or case 2), then we have, using (\ref{case1}),
\begin{eqnarray*}
\Dt &=& \sum_{j=1}^n \gamma_Q(X_j)\nabla_{X_j}^{\Si M\otimes \Si N} \\
&=& \sum_{j=1}^n \left(\gamma_M(X_j)\otimes\Id\right)
\nabla_{X_j}^{\Si M\otimes \Si N} \\
&=& \DM \\
&=& \DMt .
\end{eqnarray*}
If $n$ is odd and $m$ is even (case 3) we get from (\ref{case3}) on $\Si M
\otimes \Si^+ N$
$$
\Dt = \DM = \DMt
$$
whereas on $\Si M \otimes \Si^- N$ we obtain
$$
\Dt = -\DM = -\DMt .
$$
Finally, if $n$ and $m$ are both odd (case 4) we get from (\ref{case4})
$$
\Dt = i\left(
\begin{array}{cc}
0 & \DM \\
-\DM & 0
\end{array}
\right) .
$$
In all cases we see that $\Dt$ is formally self-adjoint because $\DM$ is,
and $\Dt^2 = (\DMt)^2$.
\qed

%%%%%%%%%%%%%%%%%%%%%%%%%%%%%%%%%%%%%%%%%%%%%%%%%%%%%%%%%%%%%%%%%%%%%%
%%%%%%%%%%%%%%%%%%%%%%%%%%%%%%%%%%%%%%%%%%%%%%%%%%%%%%%%%%%%%%%%%%%%%%

\section{The Estimate in Arbitrary Codimension}

With the preparations of the previous two sections we are now able
to bound Dirac eigenvalues of the submanifold $M$ by extrinsic
data provided $Q$ is a very ``nice'' manifold meaning that it
carries parallel or, more generally, Killing spinors.

Recall that a spinor $\psi$ on $Q$ is called {\em Killing spinor} with
{\em Killing constant} $\alpha\in\C$ if
$$
\nabla^{\Si Q}_X\psi = \alpha \cdot \gamma_Q(X) \psi
$$
for all $X \in TQ$.
Obviously, the set of Killing spinors with Killing constant $\alpha$ forms
a vector space.
We define
$$
\nu(Q,\alpha) := \dim\{\psi |\ \psi \mbox{ is a spinor on $Q$ satisfying }
\nabla^{\Si Q}_X\psi = \alpha \cdot \gamma_Q(X) \psi \} .
$$
Let $\mu(Q,n,\alpha)$ be the smallest integer greater than or equal to
$\nu(Q,\alpha)/2$ if dimension $n$ and codimension $m$
of $M$ are both odd; put $\mu(Q,n,\alpha) := \nu(Q,\alpha)$ otherwise.
\ep

\begin{theorem}\label{abschreell}
Let $Q$ be a (not necessarily complete) Riemannian spin manifold.
Let $\alpha \in \R$.
Let $M$ be an $n$-dimensional closed Riemannian spin manifold isometrically
immersed in $Q$.
Let the normal bundle $N$ of $M$ in $Q$ carry the induced spin
structure.
Let $H$ be the mean curvature vector field of $M$ in $Q$.

Then there are at least $\mu=\mu(Q,n,\alpha)$ eigenvalues $\lambda_1,\ldots,
\lambda_\mu$ of the Dirac operator $\DM$ on $M$ with coefficients
in the spinor bundle of the normal bundle, $\Si N$, satisfying the
estimate
$$
\lambda_j^2 \le n^2\alpha^2 + \frac{n^2}{4\cdot\vol(M)}\int_M |H|^2
$$
\end{theorem}
\ep

Manifolds with parallel spinors can be characterized in terms of
holonomy \cite{Hitchin,Wang}.
Simply connected manifolds with Killing spinors are described in
\cite{Baer2}.
Let us look at the most prominent examples.
\ep

{\bf Example.}
Let $Q = \R^{n+m}$ with the Euclidean metric.
The spinor bundle $\Si\R^{n+m}$ can be trivialized by parallel
spinors.
Hence for $\alpha=0$, $\nu(\R^{n+m},0)=2^{[\frac{n+m}{2}]}$, and $\mu=
\mu(\R^{n+m},n,0)=2^{[\frac{n}{2}]+[\frac{m}{2}]}$.
\ep

\begin{corollary}
Let $M$ be an $n$-dimensional closed Riemannian spin manifold,
isometrically immersed in $\R^{n+m}$.
Let the normal bundle $N$ of $M$ in $\R^{n+m}$ carry the induced
spin structure.
Let $H$ be the mean curvature vector field of $M$.

Then the Dirac operator $\DM$ of $M$ with coefficients in $\Si N$
has at least $2^{[\frac{n}{2}]+[\frac{m}{2}]}$ eigenvalues
$\lambda$ (counted with multiplicities) satisfying
$$
\lambda^2 \le \frac{n^2}{4\cdot\vol(M)} \int_M |H|^2 .
%\mbox{\qed}
$$
\qed
\end{corollary}
\ep

{\bf Example.}
Let $Q = S^{n+m}$ with the standard metric of constant sectional
curvature 1.
The spinor bundle $\Si S^{n+m}$ can be trivialized by Killing
spinors with Killing constant $\alpha=1/2$.
Hence $\nu(S^{n+m},1/2)=2^{[\frac{n+m}{2}]}$ and
$\mu=\mu(S^{n+m},n,1/2)=2^{[\frac{n}{2}]+[\frac{m}{2}]}$.
\ep

\begin{corollary}\label{qistsphaere}
Let $M$ be an $n$-dimensional closed Riemannian spin manifold,
isometrically immersed in $S^{n+m}$.
Let the normal bundle $N$ of $M$ in $S^{n+m}$ carry the induced
spin structure.
Let $H$ denote the mean curvature vector field of $M$.

Then the Dirac operator $\DM$ of $M$ with coefficients in $\Si N$
has at least $2^{[\frac{n}{2}]+[\frac{m}{2}]}$ eigenvalues
$\lambda$ satisfying
$$
\lambda^2 \le \frac{n^2}{4} + \frac{n^2}{4\cdot\vol(M)} \int_M |H|^2.
%\mbox{\qed}
$$
\qed
\end{corollary}
\ep

{\bf Proof of Theorem \ref{abschreell}.}
We will show that the Rayleigh quotient for $(\DMt)^2$ is bounded
by $n^2\alpha^2 + \frac{n^2}{4\cdot vol(M)}\int_M |H|^2$ on the
space of restrictions to $M$ of Killing spinors on $Q$ with
Killing constant $\alpha$.
Since linearly independent Killing spinors are linearly
independent at every point the space of restrictions of
Killing spinors on $Q$ to $M$
$$
\{\psi|_M\ |\ \psi \mbox{ is a spinor on $Q$ satisfying }
\nabla_X^{\Si Q}\psi = \alpha\cdot \gamma_Q(X) \psi
\hspace{0.4cm} \forall X\in TQ \}
$$
is also $\nu$-dimensional.
The factor 1/2 relating $\nu$ and $\mu$ in case $n$ and $m$ are odd comes
from the fact that in this case the operator $\DMt$ in Lemma \ref{diraclemma}
corresponds to two times the Dirac operator $\DM$.

Now let $\psi$ be a Killing spinor on $Q$ with Killing constant
$\alpha\in\R$.
Such Killing spinors have constant length and we may assume
w.l.o.g.\ that $|\psi|\equiv 1$.
We compute the Rayleigh quotient of $(\DMt)^2$ using the
notation of Lemma \ref{diraclemma}
\begin{eqnarray}
\frac{((\DMt)^2\psi,\psi)_{L^2(M)}}{(\psi,\psi)_{L^2(M)}} &=&
\frac{(\Dt^2\psi,\psi)_{L^2(M)}}{\vol(M)} \nonumber \\
&=& \frac{(\Dt\psi,\Dt\psi)_{L^2(M)}}{\vol(M)} \nonumber \\
&=& \frac{\left(\Dh\psi + \frac{n}{2}\gamma_Q(H)\psi,\Dh\psi +
\frac{n}{2}\gamma_Q(H)\psi\right)_{L^2(M)}}{\vol(M)} \nonumber \\
&=& \frac{1}{\vol(M)} \left\{ \|\Dh\psi\|^2_{\LL}
+ \frac{n}{2}(\Dh\psi,\gamma_Q(H)\psi)_{\LL}\right. \nonumber \\
&&\left. + \frac{n}{2}(\gamma_Q(H)\psi,\Dh\psi)_{\LL} +
\frac{n^2}{4}\|\gamma_Q(H)\psi\|^2_{\LL}  \right\}.
\label{rayleigh}
\end{eqnarray}
Using the Killing spinor equation we compute
\begin{eqnarray}
\Dh\psi &=& \sum_{j=1}^n \gamma_Q(X_j)\nabla_{X_j}^{\Si N}\psi
\nonumber\\
&=& \sum_{j=1}^n \gamma_Q(X_j)\alpha\gamma_Q(X_j)\psi \nonumber \\
&=& -n\cdot\alpha\cdot\psi . \label{hut}
\end{eqnarray}
Plugging (\ref{hut}) into (\ref{rayleigh}) and using skew
symmetry of $\gamma_Q(H)$ we obtain
\begin{eqnarray}
\frac{((\DMt)^2\psi,\psi)_{\LL}}{(\psi,\psi)_{\LL}} &=&
\frac{1}{\vol(M)}\left\{ n^2\alpha^2\vol(M) -
\frac{n^2\alpha}{2}(\psi,\gamma_Q(H)\psi)_{\LL} - \right.\nonumber\\
&&\left.\frac{n^2\alpha}{2}(\gamma_Q(H)\psi,\psi)_{\LL}
-\frac{n^2}{4} (\psi,\gamma_Q(H)\gamma_Q(H)\psi)_{\LL}   \right\}
\nonumber \\
&=& n^2\alpha^2 + \frac{n^2}{4}\frac{\int_M |H|^2|\psi|^2}{\vol(M)}\nonumber\\
&=& n^2\alpha^2 + \frac{n^2}{4}\frac{\int_M |H|^2}{\vol(M)} .
\end{eqnarray}
Since the Rayleigh quotient of $(\DMt)^2$ is bounded by
$n^2\alpha^2 + \frac{n^2}{4}\frac{\int_M |H|^2}{vol(M)}$ on a
$\nu$-dimensional space of spinors on $M$ the min-max
principle implies the assertion.
\qed
\ep

{\bf Remark.}
Corollary \ref{qistsphaere} can also be derived by looking at the
immersion $M \hookrightarrow S^{n+m} \subset \R^{n+m+1}$ and using
the parallel spinors on $\R^{n+m+1}$.
This is a general fact; Killing spinors with nonzero real Killing constant
on a manifold $Q$ correspond to parallel spinors on the cone over $Q$
\cite{Baer2}.
This way one can avoid dealing with Killing spinors for real Killing
constant.

There are also manifolds with Killing spinors for purely imaginary
Killing constant $\alpha$.
We can still get an eigenvalue estimate but we have to replace the
$L^2$-norm of the mean curvature $H$ by the $L^\infty$-norm.
\ep

\begin{theorem}\label{abschimaginaer}
Let $Q$ be a (not necessarily complete) Riemannian spin manifold.
Let $\alpha \in i\R$.
Let $M$ be an $n$-dimensional closed Riemannian spin manifold isometrically
immersed in $Q$.
Let the normal bundle $N$ of $M$ in $Q$ carry the induced spin
structure.
Let $H$ be the mean curvature vector field of $M$ in $Q$.

Then there are at least $\mu=\mu(Q,n,\alpha)$ eigenvalues $\lambda_1,\ldots,
\lambda_\mu$ of the Dirac operator $\DM$ on $M$ with coefficients
in the spinor bundle of the normal bundle, $\Si N$, satisfying the
estimate
$$
|\lambda_j| \le n\left( |\alpha| + \frac{1}{2}\| H\|_{L^\infty(M)} \right).
$$
\end{theorem}
\ep

{\bf Proof.}
We take a Killing spinor $\psi$ on $Q$ for the Killing constant $\alpha$
and plug $\psi|_M$ into the Rayleigh quotient of $(\DMt)^2$.
The same computations as in the proof of Theorem \ref{abschreell} yield
\begin{eqnarray}
\frac{((\DMt)^2\psi,\psi)_{\LL}}{(\psi,\psi)_{\LL}} &=&
\frac{1}{\|\psi\|^2_{\LL}}\left\{ n^2|\alpha|^2\|\psi\|^2_{\LL} -
\frac{n^2\alpha}{2}(\psi,\gamma_Q(H)\psi)_{\LL} - \right.\nonumber\\
&&\left.\frac{n^2\bar{\alpha}}{2}(\gamma_Q(H)\psi,\psi)_{\LL}
-\frac{n^2}{4} (\psi,\gamma_Q(H)\gamma_Q(H)\psi)_{\LL}   \right\}
\nonumber \\
&=& n^2|\alpha|^2 - n^2\alpha
\frac{(\psi,\gamma_Q(H)\psi)_{\LL}}{\|\psi\|^2_{\LL}} + \frac{n^2}{4}
\frac{\int_M |H|^2|\psi|^2}{\|\psi\|^2_{\LL}}
\label{rayimag}
\end{eqnarray}
Note that $|\psi|$ is no longer constant.
For any $\beta > 0$ we estimate
\begin{eqnarray}
\left| \alpha (\psi,\gamma_Q(H)\psi)_{\LL} \right| &\le&
|\alpha| \cdot \int_M |\psi|^2 |H| \nonumber\\
&=& \int_M \beta |\alpha| |\psi| \cdot \beta^{-1} |\psi| |H| \nonumber\\
&\le& \sqrt{\beta^2 \int_M |\alpha|^2 |\psi|^2} \cdot
\sqrt{\beta^{-2} \int_M |\psi|^2 |H|^2} \nonumber\\
&\le& \frac{1}{2} \left\{ \beta^2 |\alpha|^2\int_M |\psi|^2 +
\beta^{-2} \int_M |\psi|^2 |H|^2 \right\} \nonumber\\
&\le& \frac{\|\psi\|^2_{\LL}}{2}\left\{ \beta^2 |\alpha|^2 +
\beta^{-2} \|H\|^2_{L^\infty(M)} \right\} .
\label{schranke}
\end{eqnarray}
For $\beta^2=\frac{\|H\|_{L^\infty(M)}}{|\alpha|}$ inequality (\ref{schranke})
yields
\begin{equation}
\left| \alpha (\psi,\gamma_Q(H)\psi)_{\LL} \right| \le
\|\psi\|^2_{\LL}\cdot|\alpha| \cdot\|H\|_{L^\infty(M)} .
\label{schranke2}
\end{equation}
Plugging (\ref{schranke2}) into (\ref{rayimag}) gives
\begin{eqnarray*}
\frac{((\DMt)^2\psi,\psi)_{\LL}}{(\psi,\psi)_{\LL}} &\le&
n^2\cdot|\alpha|^2 + n^2\cdot|\alpha|\cdot\|H\|_{L^\infty(M)} +
\frac{n^2}{4}\cdot\|H\|_{L^\infty(M)}^2\\
&=& n^2\cdot\left( |\alpha| + \frac{1}{2} \|H\|_{L^\infty(M)}\right)^2 .
\end{eqnarray*}
The min-max principle implies the assertion.
\qed
\ep

A discussion of manifolds with Killing spinors for imaginary Killing
constant can be found in \cite{Baum2}.
The most important example is hyperbolic space $Q = H^{n+m}$.
In this case for $\alpha=\frac{i}{2}$ we have
$\nu(H^{n+m},i/2)=2^{[\frac{n+m}{2}]}$ and
$\mu=\mu(H^{n+m},n,i/2)=2^{[\frac{n}{2}]+[\frac{m}{2}]}$.
We obtain
\ep

\begin{corollary}
\label{hypkor}
Let $M$ be an $n$-dimensional closed Riemannian spin manifold,
isometrically immersed in $H^{n+m}$.
Let the normal bundle $N$ of $M$ in $H^{n+m}$ carry the induced
spin structure.
Let $H$ denote the mean curvature vector field of $M$.

Then the Dirac operator $\DM$ of $M$ with coefficients in $\Si N$
has at least $2^{[\frac{n}{2}]+[\frac{m}{2}]}$ eigenvalues
$\lambda$ satisfying
$$
|\lambda| \le \frac{n}{2}\left( 1 + \| H\|_{L^\infty(M)} \right).
%\mbox{\qed}
$$
\qed
\end{corollary}

{\bf Remark.}
If we introduce the {\em extrinsic radius} of $M$ in $Q$,
$$
rad_Q(M) = \inf \{ R>0\ |\ \exists\ p \in X \mbox{ s.t. }M \subset B(p,R) \},
$$
then the estimate in Theorem \ref{abschimaginaer} can be replaced by
$$
|\lambda_j| \le n\cdot\left( |\alpha| + \frac{e^{|\alpha|rad_Q(M)}}{2}
\sqrt{\frac{\int_M |H|^2}{\vol(M)}}\right) .
$$
Similarly, in Corollary \ref{hypkor} we also obtain the estimate
$$
|\lambda| \le \frac{n}{2}\left( 1 + e^{rad_{H^{n+m}}(M)/2}
\sqrt{\frac{\int_M |H|^2}{\vol(M)}} \right).
$$
The proof is a variation of that of Theorem \ref{abschimaginaer}
using a simple control of the growth of Killing spinors along
geodesics.
The details are left to the reader.

%%%%%%%%%%%%%%%%%%%%%%%%%%%%%%%%%%%%%%%%%%%%%%%%%%%%%%%%%%%%%%%%%%%%%%
%%%%%%%%%%%%%%%%%%%%%%%%%%%%%%%%%%%%%%%%%%%%%%%%%%%%%%%%%%%%%%%%%%%%%%

\section{Hypersurfaces}

We now turn to hypersurfaces, i.e.\ to the case $m=1$.
The point is that now the normal bundle $N$ is an oriented real
line bundle, hence trivial.
Therefore $\Si N$ is the trivial complex line bundle and
$\DM = D_M$ is the classical (untwisted) Dirac operator on $M$ acting on
spinors.
Thus Theorem \ref{abschreell} immediately gives
\ep

\begin{theorem}\label{abschreellhyp}
Let $Q$ be a (not necessarily complete) Riemannian spin manifold.
Let $\alpha \in \R$.
Let $M$ be an $n$-dimensional closed oriented hypersurface isometrically
immersed in $Q$.
Let $M$ carry the induced spin structure.
Let $H$ be the mean curvature of $M$ in $Q$.

Then there are at least $\mu=\mu(Q,n,\alpha)$ eigenvalues $\lambda_1,\ldots,
\lambda_\mu$ of the classical Dirac operator $D_M$ on $M$ satisfying the
estimate
$$
\lambda_j^2 \le n^2\alpha^2 + \frac{n^2}{4\cdot\vol(M)}\int_M H^2 .
$$
\qed
\end{theorem}
\ep

Looking at special cases for $Q$ we get corollaries corresponding to those of
the previous section.
\ep

\begin{corollary}\label{hypeukl}
Let $M$ be an $n$-dimensional closed oriented hypersurface
isometrically immersed in $\R^{n+1}$.
Let $M$ carry the induced spin structure.
Let $H$ be the mean curvature of $M$ in $\R^{n+1}$.

Then the classical Dirac operator $D_M$ of $M$
has at least $2^{[\frac{n}{2}]}$ eigenvalues
$\lambda$ (counted with multiplicities) satisfying
$$
\lambda^2 \le \frac{n^2}{4\cdot\vol(M)} \int_M H^2 .
%\mbox{\qed}
$$
\qed
\end{corollary}
\ep

Note that equality is attained for the standard sphere $M=S^n\subset\R^{n+1}$.
Indeed, $S^n$ has the eigenvalue $\lambda=n/2$ with multiplicity $2^{[n/2]}$.
It is interesting to combine the upper eigenvalue estimate with a lower
bound if $M=S^2$.
Namely, if $M$ is a surface diffeomorphic to $S^2$, then all eigenvalues
of the Dirac operator on $M$ satisfy
\begin{equation}
\lambda^2 \ge \frac{4\cdot\pi}{\mbox{area}(M)}
\nonumber
\label{meine}
\end{equation}
with equality for the eigenvalue of smallest absolute value if and only if
the metric of $M$ has constant Gauss curvature \cite[Thm.~2]{Baer3}.
Combining this with Corollary \ref{hypeukl} yields
$$
\frac{1}{\mbox{area}(M)} \int_M H^2 \ge \lambda^2 \ge \frac{4\cdot\pi}{\mbox{area}(M)}
$$
for the ``smallest'' Dirac eigenvalue $\lambda$ and hence in particular
\begin{equation}
\int_M H^2 \ge 4\cdot\pi .
\nonumber
\label{willmore}
\end{equation}

This is known as the {\em Willmore inequality} \cite[Thm.~7.2.2]{Willmore}.
It actually holds for all closed oriented surfaces immersed in $\R^3$.
The equality case in (\ref{willmore}) can also easily be discussed.
%If $\int_M H^2 = 4\cdot\pi$, then equality must hold in (\ref{meine}) and
%hence $M$ must be a sphere of constant Gauss curvature.
%After a rescaling we may assume that the Gauss curvature is~1.
%Eigenspinors $\phi$ on $S^2$ for the eigenvalue $\lambda=n/2=1$ are known to
%be Killing spinors
%\begin{equation}
%\nabla^{S^2}_X\phi = \pm \frac{1}{2}\gamma_{S^2}(X)\phi .
%\nonumber
%\label{willkill}
%\end{equation}

%On the other hand, from the proof of Corollary \ref{hypeukl} we know that
%they are restrictions of parallel spinors on $\R^3$ and hence by
%(\ref{gaussformel6})
%\begin{equation}
%0 = \nabla^{\R^3}_X \phi = \nabla^{S^2}_X \phi + \frac{1}{2}
%\gamma_{S^2}(B(X))\phi
%\nonumber
%\label{willparallel}
%\end{equation}
%for all tangent vectors $X$.
%From (\ref{willkill}) and (\ref{willparallel}) we get
%$$
%B = \pm\Id .
%$$
%Here $B$ is the shape operator of $M\subset\R^3$.
%Hence we have shown that the first and the second fundamental form of $M$
%are the same (up to rescaling) as the ones of the standard sphere.
%The fundamental theorem of surface theory now implies that $M$ is the
%standard sphere in $\R^3$ up to congruence and rescaling.

There is a conjecture that if $M$ is an embedded torus, then
$$
\int_M H^2 \ge 2\cdot\pi^2 .
$$
It is established for some classes of tori
\cite[Thm.~7.2.4]{Willmore}, \cite{ST}, \cite{HJP}, \cite{LS}, \cite{LY}
but in full generality it is still open.
\ep

{\bf Question.}
Can one show that Dirac eigenvalues $\lambda$ of tori isometrically embedded in
$\R^3$ satisfy
$$
\lambda^2 \ge \frac{2\cdot\pi^2}{\mbox{area}(M)}
$$
or at least
$$
\lambda^2 \ge \frac{4\cdot\pi}{\mbox{area}(M)}
\hspace{0.3cm}?
$$
\ep

It is known that the spin structure on a 2-torus induced by an {\em embedding}
in $\R^3$ is never trivial and thus $\lambda\not= 0$.
In contrast, the spin structure on a torus induced by an {\em immersion}
can be trivial and hence $\lambda=0$ can occur.
Therefore the answer to the question is certainly ``no'' for immersed tori.
\ep

\begin{corollary}\label{hypsphaere}
Let $M$ be an $n$-dimensional closed oriented hypersurface
isometrically immersed in $S^{n+1}$.
Let $M$ carry the induced spin structure.
Let $H$ denote the mean curvature of $M$ in $S^{n+1}$.

Then the classical Dirac operator $D_M$ of $M$
has at least $2^{[\frac{n}{2}]}$ eigenvalues
$\lambda$ satisfying
$$
\lambda^2 \le \frac{n^2}{4} + \frac{n^2}{4\cdot\vol(M)} \int_M H^2.
%\mbox{\qed}
$$
\qed
\end{corollary}
\ep

Equality in Corollary \ref{hypsphaere} is attained for hyperspheres
in $S^{n+1}$ cut out by affine hyperplanes of $\R^{n+2}$.
For example, an equatorial hypersphere $S^n \subset S^{n+1}$ is
totally geodesic, hence $H \equiv 0$, and $\lambda=\frac{n}{2}$
is an eigenvalue.

By combining Corollary \ref{hypsphaere} with (\ref{meine}) we get
a Willmore inequality for 2-spheres immersed in $S^3$.
%\begin{corollary}\label{willmoresphere}
Namely, let $M$ be a closed oriented surface of genus 0
isometrically immersed in $S^{3}$.
Let $H$ denote the mean curvature of $M$ in $S^{3}$.
Then
$$
4\cdot\pi \le \mbox{area}(M) + \int_M H^2 .
%\mbox{\qed}
$$
%\qed
%\end{corollary}
%\ep

Theorem \ref{abschimaginaer} yields in the case of hypersurfaces
\ep

\begin{theorem}\label{abschimaginaerhyp}
Let $Q$ be a (not necessarily complete) Riemannian spin manifold.
Let $\alpha \in i\R$.
Let $M$ be an $n$-dimensional closed oriented hypersurface isometrically
immersed in $Q$.
Let $M$ carry the induced spin structure.
Let $H$ be the mean curvature of $M$ in $Q$.

Then there are at least $\mu=\mu(Q,n,\alpha)$ eigenvalues $\lambda_1,\ldots,
\lambda_\mu$ of the classical Dirac operator $D_M$ on $M$ satisfying the
estimate
$$
|\lambda_j| \le n\left( |\alpha| + \frac{1}{2}\| H\|_{L^\infty(M)} \right).
$$
\qed
\end{theorem}
\ep

\begin{corollary}\label{hyphyp}
Let $M$ be an $n$-dimensional closed oriented hypersurface
isometrically immersed in $H^{n+1}$.
Let $M$ carry the induced spin structure.
Let $H$ denote the mean curvature of $M$ in $H^{n+1}$.

Then the classical Dirac operator $D_M$ of $M$
has at least $2^{[\frac{n}{2}]}$ eigenvalues
$\lambda$ satisfying
$$
|\lambda| \le \frac{n}{2}\left( 1 + \| H\|_{L^\infty(M)} \right).
%\mbox{\qed}
$$
\qed
\end{corollary}

If $M$ is the distance sphere from a point $p\in H^{n+1}$ with radius
$r>0$, then $H = \coth(r)$ and the smallest positive eigenvalue of the
Dirac operator is $\lambda(r) = \frac{n}{2\cdot\sinh(r)}$.
Hence the estimate in Corollary \ref{hyphyp} is asymptotically sharp for
$r \to 0$ in the sense that
$$
\lim_{r\searrow 0} \frac{|\lambda(r)|}{\frac{n}{2}\left( 1 + \| H\|_{L^\infty(M)}
\right)}
= 1 .
$$
Of course, we also obtain a version of the Willmore inequality.
%\begin{corollary}\label{willmorehyp}
Namely, let $M$ be a closed oriented surface of genus 0
isometrically immersed in $H^{3}$.
Let $H$ denote the mean curvature of $M$ in $H^{3}$.
Then
$$
4\cdot\pi \le \left( 1 + \| H\|_{L^\infty(M)} \right)^2\cdot \mbox{area}(M) .
%\mbox{\qed}
$$
%\qed
%\end{corollary}
%\ep

As mentioned at the end of Section 3 the $L^\infty$-norm of $H$ in the
estimates of Theorem \ref{abschimaginaerhyp} and Corollary \ref{hyphyp} can
be replaced by the $L^2$-norm if we insert an additional term involving
the extrinsic radius of $M$ in the ambient space.

Before concluding this section let us note that Theorem
\ref{abschreellhyp} can be improved if $M$ bounds a relatively compact
open subset $W\subset Q$.
To do this we need the following version of the standard variational
characterization of eigenvalues.
\ep

\begin{lemma}\label{varivers}
Let $\LLL$ be a separable Hilbert space, let $\HHH_1$ and $\HHH_2$ be two
$\nu$-dimen\-sional subspaces, orthogonal to each other.
Let $A$ be a nonnegative selfadjoint operator on $\LLL$ with pure point spectrum
$0 \le \lambda_1 \le \lambda_2 \le \ldots$ where the eigenvalues are repeated
according to their multiplicity.

Denote the Rayleigh quotient of $A$ by $Q^A$:
$$
Q^A(\phi) = \frac{(A\phi,\phi)}{\|\phi\|^2}.
$$
Let $Q^A$ be bounded on $\HHH_1$ and $\HHH_2$ by some constant $C$.

Then for $j=1,\ldots,\nu$:
$$
\frac{\lambda_j+\lambda_{2\nu-j+1}}{2} \le C .
$$
\end{lemma}
\ep

Note that this result is sharp as one sees from the following example.
Let $\LLL = \C^4$ with the standard orthonormal basis $e_1,e_2,e_3,e_4$.
Let
\[
A = \left(
\begin{array}{cccc}
1&0&0&0\\
0&2&0&0\\
0&0&3&0\\
0&0&0&4
\end{array}
\right) .
\]
Let $\HHH_1$ be spanned by $e_2 + e_3$ and $e_1+e_4$ and $\HHH_2$ by
$e_2 - e_3$ and $-e_1+e_4$.
One checks that $Q^A$ is bounded by $C = 5/2$ on $\HHH_i$.

Indeed, $\frac{\lambda_1+\lambda_4}{2} = C$ and
$\frac{\lambda_2+\lambda_3}{2} = C$.
Note that the estimate already fails for $\frac{\lambda_2+\lambda_4}{2}$.
\ep

{\bf Proof of Lemma \ref{varivers}.}
Let $\phi_1,\phi_2,\ldots$ be an orthonormal basis of $\LLL$ consisting
of eigenvectors of $A$ for the eigenvalues $\lambda_1,\lambda_2,\ldots$ .
For $j\le k$ denote the span of $\phi_j,\phi_{j+1},\ldots,\phi_k$ by $E_j^k$
and the span of $\phi_j,\phi_{j+1},\ldots$ by $E_j^\infty$.

Let $1\le j \le \nu$.
>From $\dim(\HHH_i)=\nu$ and $\codim(E_j^\infty)=j-1$ we see that
$\dim(\HHH_i\cap E_j^\infty)\ge \nu - j + 1$ and hence
$$
\dim\left((\HHH_1\cap E_j^\infty)\oplus(\HHH_2\cap E_j^\infty)\right)
\ge 2\nu - 2j + 2 .
$$
We look at the map
$$
\Phi : (\HHH_1\cap E_j^\infty)\oplus(\HHH_2\cap E_j^\infty)
\rightarrow E_j^{2\nu -j} ,
$$
$$
\Phi(\psi_1\oplus\psi_2)= (\mbox{orthogonal projection onto }E_j^{2\nu -j})
(\psi_1-\psi_2) .
$$
Since $\dim(E_j^{2\nu -j})=2\nu - 2j +1$ the kernel of $\Phi$ must be nontrivial.
Let $\psi_1\oplus\psi_2$ be in this kernel, $\psi_i\in\HHH_i\cap E_j^\infty$.
We express $\psi_1$ and $\psi_2$ in the basis $\phi_1,\phi_2,\ldots$,
$$
\psi_i = \sum_{k=j}^\infty \alpha_{i,k}\phi_k .
$$
Since $\psi_1\oplus\psi_2$ is in the kernel of $\Phi$ we have $\alpha_{1,k}
=\alpha_{2,k}$ for $k=j,\ldots,2\nu-j$.
Hence
\begin{eqnarray}
Q^A(\psi_1+\psi_2) &=& \frac{(A(\psi_1+\psi_2),\psi_1+\psi_2)}
{\|\psi_1+\psi_2\|^2}
\nonumber \\
&=&\frac{(\sum_{k\ge j}(\alpha_{1,k}+\alpha_{2,k})\lambda_k\phi_k,
\sum_{l\ge j}(\alpha_{1,l}+\alpha_{2,l})\phi_l)}
{(\sum_{k\ge j}(\alpha_{1,k}+\alpha_{2,k})\phi_k,
\sum_{l\ge j}(\alpha_{1,l}+\alpha_{2,l})\phi_l)}
\nonumber \\
&=& \frac{\sum_{k\ge j}|\alpha_{1,k}+\alpha_{2,k}|^2\lambda_k}
{\sum_{k\ge j}|\alpha_{1,k}+\alpha_{2,k}|^2}
\nonumber \\
&\ge& \frac{\lambda_j\sum_{k\ge j}|\alpha_{1,k}+\alpha_{2,k}|^2}
{\sum_{k\ge j}|\alpha_{1,k}+\alpha_{2,k}|^2}
\nonumber \\
&=& \lambda_j
\label{qplus}
\end{eqnarray}
and similarly
\begin{equation}
Q^A(\psi_1-\psi_2) =
\frac{\sum_{k\ge 2\nu-j+1}|\alpha_{1,k}-\alpha_{2,k}|^2\lambda_k}
{\sum_{k\ge 2\nu-j+1}|\alpha_{1,k}-\alpha_{2,k}|^2}
\ge \lambda_{2\nu-j+1} .
\label{qminus}
\end{equation}
Adding (\ref{qplus}) and (\ref{qminus}) we get, using the fact that
$\psi_1$ and $\psi_2$ are orthogonal,
\begin{eqnarray}
\lambda_j + \lambda_{2\nu-j+1} &\le&
Q^A(\psi_1+\psi_2) + Q^A(\psi_1-\psi_2)
\nonumber \\
&=& \frac{(A\psi_1,\psi_1) + (A\psi_1,\psi_2) +
(A\psi_2,\psi_1) + (A\psi_2,\psi_2)}
{\|\psi_1\|^2 + \|\psi_2\|^2}
\nonumber \\
&&
+ \frac{(A\psi_1,\psi_1) - (A\psi_1,\psi_2) -
(A\psi_2,\psi_1) + (A\psi_2,\psi_2)}
{\|\psi_1\|^2 + \|\psi_2\|^2}
\nonumber \\
&=& 2\cdot\frac{(A\psi_1,\psi_1) + (A\psi_2,\psi_2)}
{\|\psi_1\|^2 + \|\psi_2\|^2}
\nonumber \\
&\le& 2\cdot\frac{C\cdot\|\psi_1\|^2 + C\cdot\|\psi_2\|^2}
{\|\psi_1\|^2 + \|\psi_2\|^2}
\nonumber \\
&=& 2\cdot C .
\nonumber
\end{eqnarray}
\qed
\ep

\begin{theorem}\label{hoehere}
Let $Q$ be a (not necessarily complete) $(n+1)$-dimensional Riemannian
spin manifold.
Let $\alpha \in \R$.
Let $W\subset Q$ be a relatively compact open subset with smooth boundary $M$.
Let $M$ carry the induced spin structure.
Let $H$ be the mean curvature of $M$ in $Q$.

Then there are at least $2\mu=2\mu(Q,n,\alpha)$ eigenvalues $\lambda_1,
\ldots,\lambda_{2\mu}$ of the classical Dirac operator $D_M$ on $M$
satisfying the estimate
$$
\frac{\lambda_j^2+\lambda_{2\mu-j+1}^2}{2} \le n^2\alpha^2 +
\frac{n^2}{4\cdot\vol(M)} \int_M H^2
$$
for $j=1,\ldots,\mu$.
\end{theorem}
\ep

{\bf Remark.}
Theorem \ref{hoehere} should be read as follows.
By Theorem \ref{abschreellhyp} one knows that there are $\mu$
eigenvalues $\lambda_1,\ldots,\lambda_\mu$ satisfying
$$
\lambda_j^2 \le n^2\alpha^2 + \frac{n^2}{4\cdot\vol(M)}\int_M H^2
=: C .
$$
Now Theorem \ref{hoehere} says in particular that
$$
\lambda_{\mu+j}^2 \le 2C
$$
for $1\le j \le \mu$.
If one of the ``higher'' eigenvalues $\lambda_{\mu+j}^2$ is much larger
that $C$, say $\lambda_{\mu+j}^2 \ge C+\epsilon$, then the
corresponding ``small'' eigenvalue $\lambda_{\mu-j+1}^2$ must satisfy
the stronger estimate
$$
\lambda_{\mu-j+1}^2 \le C-\epsilon .
$$

{\bf Example.}
If $M=S^n\subset Q=\R^{n+1}$ is the standard sphere, then
by Corollary~\ref{hypeukl} there are at least $\mu=2^{[n/2]}$ eigenvalues
$\lambda_1,\ldots,\lambda_\mu$ satisfying
$$
\lambda_j^2 \le \frac{n^2}{4\cdot\vol(M)}\int_M H^2 = \frac{n^2}{4} .
$$
The scalar curvature of $M$ is $n(n-1)$ and hence by Friedrich's estimate
\cite{Friedrich} all eigenvalues $\lambda$ satisfy
$$
\lambda^2 \ge \frac{1}{4}\frac{n}{n-1}n(n-1) = \frac{n^2}{4} .
$$
Thus by Theorem \ref{hoehere} there must actually be $2\mu=2^{[n/2]+1}$
eigenvalues $\lambda_j$ satisfying
$$
\lambda_j^2 \le \frac{n^2}{4} .
$$
In fact, the eigenvalues $\lambda=\pm\frac{n}{2}$ both have multiplicity
$2^{[n/2]}$.
\ep

{\bf Proof.}
We know that the Rayleigh quotient of $A=(\Dt_M)^2=(\DMt)^2$ is bounded by
$$
C=n^2\alpha^2 + \frac{n^2}{4\cdot\vol(M)}\int_M H^2
$$
on the $\nu$-dimensional subspace
$$
{\mbox{\goth H}}_1 =
\{\psi|_M\ |\ \psi \mbox{ is a spinor on $Q$ satisfying }
\nabla_X^{\Si Q}\psi = \alpha\cdot \gamma_Q(X) \psi
\hspace{0.4cm} \forall X\in TQ \}
$$
of $\LLL=L^2(M,\Si Q|_M)$.
Let $Z$ denote the exterior unit normal field of $W \subset Q$ along
$\partial W = M$.
Define another subspace of $L^2(M,\Si Q|_M)$ by
$$
{\mbox{\goth H}}_2 :=
\{ \gamma_Q(Z)\psi\ |\ \psi \in {\mbox{\goth H}}_1 \}.
$$
Since Clifford multiplication with $Z$ is an isomorphism (with
inverse $-\gamma_Q(Z)$) ${\mbox{\goth H}}_2$ is also $\nu$-dimensional.

We compute the Rayleigh quotient on ${\mbox{\goth H}}_2$.
Let $\gamma_Q(Z)\psi\in \mbox{\goth H}_2$.
\begin{eqnarray*}
\frac{\left(\Dt_M(\gamma_Q(Z)\psi),\Dt_M(\gamma_Q(Z)\psi)\right)_{L^2(M)}}
{\left(\gamma_Q(Z)\psi,\gamma_Q(Z)\psi\right)_{L^2(M)}} &=&
\frac{\left(-\gamma_Q(Z)\Dt_M\psi,-\gamma_Q(Z)\Dt_M\psi\right)_{L^2(M)}}
{(\psi,\psi)_{L^2(M)}} \\
&=& \frac{(\Dt_M\psi,\Dt_M\psi)_{L^2(M)}}
{(\psi,\psi)_{L^2(M)}} \\
&=& n^2\alpha^2 + \frac{n^2}{4\cdot\vol(M)}\int_M H^2 .
\end{eqnarray*}
It remains to show that ${\mbox{\goth H}}_1$ and ${\mbox{\goth H}}_2$ are
$L^2$-orthogonal.
This is where we use that $M$ bounds.
Let $\psi_1$ and $\psi_2$ be in ${\mbox{\goth H}}_1$.
Using the Green's formula \cite[p.~115, eq.~(5.7)]{LM} we get
\begin{eqnarray*}
\int_M \langle \psi_1,\gamma_Q(Z)\psi_2 \rangle &=&
\left(\psi_1,D_W\psi_2\right)_{L^2(W)} -
\left(D_W\psi_1,\psi_2\right)_{L^2(W)} \\
&=& \left(\psi_1,-(n+1)\alpha\psi_2\right)_{L^2(W)} -
\left(-(n+1)\alpha\psi_1,\psi_2\right)_{L^2(W)} \\
&=& 0 .
\end{eqnarray*}
Lemma \ref{varivers} yields the assertion.
\qed
\ep

Note that Theorem \ref{hoehere} fails if the hypersurface $M$ does
not bound in $Q$.
For example, for the flat 3-torus $Q = \R^3/\Z^3$ we have
$\nu(\R^3/\Z^3,0)=2$ and thus $\mu(\R^3/\Z^3,2,0)=2$.
Let $M=\R^2/\Z^2 \subset Q=\R^3/\Z^3$ be a linear subtorus.
Then $M$ is totally geodesic, hence $H\equiv 0$.
Theorem \ref{abschreellhyp} says that the eigenvalue $\lambda=0$
has multiplicity 2 at least.
Indeed, the multiplicity is precisely 2.
If Theorem \ref{hoehere} could be applied in this case it would say
that the multiplicity of $\lambda=0$ is at least 4 which is not true.

%%%%%%%%%%%%%%%%%%%%%%%%%%%%%%%%%%%%%%%%%%%%%%%%%%%%%%%%%%%%%%%%%%%%%%
%%%%%%%%%%%%%%%%%%%%%%%%%%%%%%%%%%%%%%%%%%%%%%%%%%%%%%%%%%%%%%%%%%%%%%

\section{Higher Eigenvalues}

So far we have given estimates on the smallest $\mu$ eigenvalues only.
It is also possible to bound higher eigenvalues.
We will show how to obtain bounds on higher Dirac eigenvalues which
involve the eigenvalues of the Laplace-Beltrami operator acting on
functions.

We return to the case of arbitrary codimension.
\ep

\begin{theorem}
Let $Q$ be a Riemannian spin manifold.
Let $\alpha \in \R$.
Let $M$ be an $n$-dimensional closed Riemannian spin manifold isometrically
immersed in $Q$.
Let the normal bundle $N$ of $M$ in $Q$ carry the induced spin
structure.
Let $H$ be the mean curvature vector field of $M$ in $Q$.
Write $\mu=\mu(Q,n,\alpha)$.
Denote the eigenvalues of the Laplace-Beltrami operator $\Delta$
acting on functions on $M$ by
$$
0 = \lambda_0(\Delta) < \lambda_1(\Delta) \le \lambda_2(\Delta)
\le \lambda_3(\Delta) \le \ldots
$$

Then the eigenvalues of the Dirac operator $\DM$ on $M$ with coefficients
in the spinor bundle of the normal bundle, $\Si N$, satisfy the estimate
$$
\lambda_{(k+1)\mu}^2 \le
\frac{n^2}{4} \| H\|^2_{L^\infty(M)} + n^2\alpha^2 + \lambda_k(\Delta)
$$
for $k=0,1,\ldots$\ .
\end{theorem}
\ep

{\bf Proof.}
Let $\psi$ be a Killing spinor on $Q$ with Killing constant $\alpha$.
Normalize $\psi$ such that $|\psi|\equiv 1$.
For a smooth function $f : M \to \R$ a computation similar to the one
in the proof of Theorem \ref{abschreell} yields
\begin{eqnarray*}
\left((\DMt)^2(f\psi),f\psi\right)_{L^2(M)} &=&
\left(\Dt(f\psi),\Dt(f\psi)\right)_{L^2(M)} \\
&=& n^2\alpha^2 (f\psi,f\psi)_{L^2(M)} \\
&& + \left( \ga_Q\left(\frac{n}{2}fH + \grad f\right)\psi,
\ga_Q\left(\frac{n}{2}fH + \grad f\right)\psi\right)_{L^2(M)} \\
&=& n^2\alpha^2 (f,f)_{L^2(M)} + \frac{n^2}{4} \int_M f^2 |H|^2
+ \left(\grad f,\grad f\right)_{L^2(M)} .
\end{eqnarray*}
For the Rayleigh quotient we obtain
\begin{eqnarray*}
\frac{\left((\DMt)^2(f\psi),f\psi\right)_{L^2(M)}}{(f\psi,f\psi)_{L^2(M)}}
&=& n^2\alpha^2 + \frac{\frac{n^2}{4} \int_M f^2 |H|^2}{(f,f)_{L^2(M)}} +
\frac{\left(\grad f,\grad f\right)_{L^2(M)}}{(f,f)_{L^2(M)}} \\
&\le& n^2\alpha^2 + \frac{n^2}{4} \| H\|^2_{L^\infty(M)} +
\frac{\left(\Delta f,f\right)_{L^2(M)}}{(f,f)_{L^2(M)}} .
\end{eqnarray*}
Using the test space spanned by products $f\psi$ where $\psi$ is a Killing
spinor on $Q$ with Killing constant $\alpha$ and $f$ is an eigenfunction
of $\Delta$ for the eigenvalue $\lambda_j(\Delta)$, $j\le k$, the min-max
principle yields the proof.
\qed

%%%%%%%%%%%%%%%%%%%%%%%%%%%%%%%%%%%%%%%%%%%%%%%%%%%%%%%%%%%%%%%%%%%%%%
%%%%%%%%%%%%%%%%%%%%%%%%%%%%%%%%%%%%%%%%%%%%%%%%%%%%%%%%%%%%%%%%%%%%%%

\ep
\ep

\noindent
Author's address:

\ep

\noindent
Mathematisches Institut\\
Universit\"at Freiburg\\
Eckerstr. 1\\
79104 Freiburg\\
Germany

\ep

\noindent
e-mail: {\tt baer@mathematik.uni-freiburg.de}

\end{document}